\documentclass{aptpub}

\authornames{S. CONNOR AND S. JACKA } 
\shorttitle{Optimal co-adapted coupling on the hypercube} 

\newcommand{\abs}[1]{\left\vert#1\right\vert}
\newcommand{\set}[1]{\left\{#1\right\}}

\newcommand{\squ}[1]{\left[#1\right]}
\newcommand{\bra}[1]{\left(#1\right)}
\newcommand{\norm}[1]{\left\Vert#1\right\Vert}
\renewcommand{\prob}[1]{\operatorname{\mathbb{P}}\left[#1\right]}

\newcommand{\Ex}[1]{\operatorname{\mathbb{E}}\left[#1\right]}

\newcommand{\indev}[1]{\mathbf{1}_{\squ{#1}}}
\newcommand{\ind}[1]{\mathbf{1}_{#1}}
\newcommand{\hyp}{\ensuremath{\mathbb{Z}_2^n}}
\newcommand{\law}{\ensuremath{\mathcal{L}}}
\newcommand{\RR}{\mathbb R}

\numberwithin{equation}{section}  

\begin{document}

\title{Optimal co-adapted coupling for the symmetric random walk on the hypercube} 

\authorone[University of Warwick]{Stephen Connor} 

\addressone{Department of Statistics, University of Warwick, Coventry, CV4 7AL. UK} 

\authorone[University of Warwick]{Saul Jacka}

\begin{abstract}
Let $X$ and $Y$ be two simple symmetric continuous-time random walks on the vertices of the $n$-dimensional hypercube, \hyp. We consider the class of co-adapted couplings of these processes, and describe an intuitive coupling which is shown to be the fastest in this class.
\end{abstract}



\section{Introduction}\label{sec:intro}

Let \hyp\ be the group of binary $n$-tuples under coordinate-wise
addition modulo 2: this can be viewed as the set of vertices of an $n$-dimensional
hypercube. For $x\in\hyp$, we write $x=\bra{x(1),\dots,x(n)}$, and define elements $\set{e_i}_0^n$ by
    \[ e_0=\bra{0,\dots,0}; \quad e_i(k) = \indev{i=k}, \,\, i=1,\dots,n \,, \]
where $\ind{}$ denotes the indicator function.
For $x,y\in\hyp$ let
    \[ \abs{x-y} = \sum_{i=1}^n\abs{x(i)-y(i)} \]
denote the Hamming distance between $x$ and $y$.

A continuous-time random walk $X$ on \hyp\ may be defined using a marked Poisson process $\Lambda$ of rate $n$, with marks distributed uniformly on the set $\set{1,2,\dots,n}$: the $i^{\text{th}}$ coordinate of $X$ is flipped  to its opposite value (zero or one) at incident times of $\Lambda$ for which the corresponding mark is equal to $i$. We write $\law\bra{X_t}$ for the law of $X$ at time $t$. The unique equilibrium distribution of $X$ is the uniform distribution on \hyp.

Suppose that we now wish to couple two such random walks, $X$ and $Y$, starting from different states.
    \begin{defn}\label{defn:coupling}
        A \emph{coupling} of $X$ and $Y$ is a process $(X',Y')$ on $\hyp\times\hyp$ such that
            \[ X'\stackrel{\mathcal{D}}{=} X \quad \text{and} \quad Y'\stackrel{\mathcal{D}}{=} Y \,. \]
        That is, viewed marginally, $X'$ behaves as a version of $X$, and $Y'$ as a version of $Y$.
    \end{defn}
For any coupling strategy $c$, write $(X_t^c,Y_t^c)$ for the value at $t$ of the pair of processes $X^c$ and $Y^c$ driven by strategy $c$, although this superscript notation may be dropped when no confusion can arise. (We assume throughout that $(X^c, Y^c)$ \emph{is} a coupling of $X$ and $Y$.) We then define the coupling time by
    \[ \tau^c = \inf\set{t\geq 0 \,:\, X^c_s = Y^c_s \;\; \forall\, s\geq t}\,. \]
Note that in general this is not necessarily a stopping time for either of the marginal processes, nor even for the joint process. For $t\geq 0$, let
    \[ U^c_t = \set{1\leq i \leq n \,:\, X^c_t(i) \neq Y^c_t(i)} \]
denote the set of unmatched coordinates at time $t$, and let
    \[ M^c_t = \set{1\leq i \leq n \,:\, X^c_t(i) = Y^c_t(i)} \]
be its complement. A simple coupling technique appears in~\cite{Aldous-1983}, and may be described as follows:
    \begin{itemize}
        \item if $X(i)$ flips at time $t$, with $i\in M_t$, then also flip coordinate $Y(i)$ at time $t$ (matched coordinates are always made to move synchronously);
        \item if $\abs{U_t}>1$ and $X(i)$ flips at time $t$, with $i\in U_t$, also flip coordinate $Y(j)$ at time $t$, where $j$ is chosen uniformly at random from the set $U_t\backslash\set{i}$;
        \item else, if $U_t=\set{i}$ contains only one element, allow coordinates $X(i)$ and $Y(i)$ to evolve independently of each other until this final match is made.
    \end{itemize}
This defines a valid coupling of $X$ and $Y$, for which existing coordinate matches are maintained and new matches made in pairs when $\abs{U_t}\geq 2$. It is also an example of a \emph{co-adapted} coupling.
    \begin{defn}\label{defn:co-adapted}
        A coupling $(X^c,Y^c)$ is called \emph{co-adapted} if there exists a filtration $\bra{\mathcal{F}_t}_{t\geq 0}$ such that
        \begin{enumerate}
            \item $X^c$ and $Y^c$ are both adapted to $\bra{\mathcal{F}_t}_{t\geq 0}$
            \item for any $0\leq s\leq t$,
                \[ \law\bra{X^c_t \,|\, \mathcal{F}_s} = \law\bra{X^c_t \,|\, X^c_s} \quad\text{and}\quad
                    \law\bra{Y^c_t \,|\, \mathcal{F}_s} = \law\bra{Y^c_t \,|\, Y^c_s}\,.
                \]
        \end{enumerate}
    \end{defn}
In other words, $(X^c,Y^c)$ is co-adapted if $X^c$ and $Y^c$ are both Markov with respect to a common filtration, $\bra{\mathcal{F}_t}_{t\geq 0}$. Note that this definition does \emph{not} imply that the joint process $(X^c,Y^c)$ is Markovian, however. If $(X^c,Y^c)$ is co-adapted then the coupling time is a randomised stopping time with respect to the individual chains, and it suffices to study the first \emph{collision} time of the two chains (since it is then always possible to make $X^c$ and $Y^c$ agree from this time onwards).


In this paper we search for the best possible coupling of the random walks $X$ and $Y$ on \hyp\ within the class $\mathcal{C}$ of all co-adapted couplings.


\section{Co-adapted couplings for random walks on \hyp}\label{sec:co-adapted}

In order to find the optimal co-adapted coupling of $X$ and $Y$, it is first necessary to be able to describe a general coupling strategy $c\in\mathcal{C}$. To this end, let $\Lambda_{ij}$ ($0\leq i,j \leq n$) be independent unit-rate
marked Poisson processes, with marks $W_{ij}$ chosen uniformly on the interval $[0,1]$. We let $\bra{\mathcal{F}_t}_{t\geq 0}$ be any filtration satisfying
    \[ \sigma\set{\,\bigcup_{i,j}\Lambda_{ij}(s),\, \bigcup_{i,j}
        W_{ij}(s)\,:\, s\leq t} \subseteq \mathcal{F}_t, \quad\forall\,t\geq 0 \,. \]
The transitions of $X^c$ and $Y^c$ will be driven by the marked Poisson processes, and controlled by a process
$\set{Q^c(t)}_{t\geq 0}$ which is adapted to $\bra{\mathcal{F}_t}_{t\geq 0}$. Here, $Q^c(t)=\set{q^c_{ij}(t)\,:\, 1\leq i,j,\leq n}$ is a $n\times n$ doubly sub-stochastic matrix. Such a matrix implicitly
defines terms $\set{q^c_{0j}(t)\,:\, 1\leq j\leq n}$ and
$\set{q^c_{i0}(t)\,:\, 1\leq i\leq n}$ such that %
    \begin{align}
        \sum_{i=0}^n q^c_{ij}(t)& =1 \quad\text{for all $1\leq j\leq n$ and $t\geq 0$\,,} \label{eqn:Q-reln1} \\
        \text{and}\quad \sum_{j=0}^n q^c_{ij}(t)& =1 \quad\text{for all $1\leq
        i\leq n$ and $t\geq 0$\,.} \label{eqn:Q-reln2}
    \end{align}
For convenience we also define $q^c_{00}(t)=0$ for all $t\geq 0$.

Note that any co-adapted coupling $(X^c,Y^c)$ must satisfy the following three constraints, all of which are due to the marginal processes $X^c(i)$ ($i=1,\dots,n$) being independent unit rate Poisson processes (and similarly for the processes $Y^c(i)$):
    \begin{enumerate}
        \item At any instant the number of jumps by the process $(X^c,Y^c)$ cannot exceed two (one on $X^c$ and one on $Y^c$); 
        \item All single and double jumps must have rates bounded above by one;
        \item For all $i=1,\dots,n$, the \emph{total} rate at which $X^c(i)$ jumps must equal one.
    \end{enumerate}

A general co-adapted coupling for $X$ and $Y$ may therefore be defined as
follows: if there is a jump in the process $\Lambda_{ij}$ at
        time $t\geq 0$, \emph{and} the mark $W_{ij}(t)$ satisfies $W_{ij}(t)\leq q_{ij}(t)$, then
        set $X^c_t=X^c_{t-}+e_i$ (mod 2) and $Y^c_t=Y^c_{t-}+e_j$ (mod 2).
Note that if $i$ (respectively $j$) equals zero, then $X^c_t=X^c_{t-}$ (respectively, $Y^c_t=Y^c_{t-}$), since $e_0=(0,\dots,0)$.

From this construction it follows directly that $X^c$ and $Y^c$ both
have the correct marginal transition rates to be continuous-time
simple random walks on \hyp\ as described above, and are co-adapted. 


\section{Optimal coupling}\label{sec:optimal}

Our proposed optimal coupling strategy, $\hat{c}$, is very simple to describe, and depends only upon the number of unmatched coordinates of $X$ and $Y$. Let $N_t = \abs{U_t}$ denote the value of this number at time $t$. Strategy $\hat{c}$ may be summarised as follows:
    \begin{itemize}
        \item matched coordinates are always made to move synchronously (thus $N^{\hat{c}}$ is a decreasing process);
        \item if $N$ is odd, all unmatched coordinates of $X$ and $Y$ are made to evolve independently until $N$ becomes even;
        \item if $N$ is even, unmatched coordinates are coupled in pairs - when an unmatched coordinate on $X$ flips (thereby making a new match), a different, uniformly chosen, unmatched coordinate on $Y$ is forced to flip at the same instant (making a total of two new matches).
    \end{itemize}
Note the similarity between $\hat{c}$ and the coupling of Aldous described in Section~\ref{sec:intro}: if $N$ is even these strategies are identical; if $N$ is odd however, $\hat{c}$ seeks to restore the parity of $N$ as fast as possible, whereas Aldous's coupling continues to couple unmatched coordinates in pairs until $N=1$.

    \begin{defn}\label{defn:proposed-optimum}
        The matrix process $\hat{Q}$ corresponding to the coupling $\hat{c}$ is as follows:
            \begin{itemize}
                \item $\hat{q}_{ii}(t)=1$ for all $i\in M_t$ and for all $t\geq 0$; \label{defn:Q-hat}
                \item if $N_t$ is odd, $\hat{q}_{i0}(t)=\hat{q}_{0i}(t) = 1$ for all $i\in U_t$;
                \item if $N_t$ is even, $\hat{q}_{i0}(t)=\hat{q}_{0i}(t) =\hat{q}_{ii}(t)= 0$ for all $i\in U_t$, and \[ \hat{q}_{ij} = \frac{1}{\abs{U_t}-1} \quad\text{for all distinct $i,j\in U_t$} \,. \]
            \end{itemize}
    \end{defn}
The coupling time under $\hat{c}$, when $(X_0,Y_0)=(x,y)$, can thus be expressed as follows:
    \begin{equation}
        \hat{\tau} = \tau^{\hat{c}} =
            \begin{cases}\label{eqn:tau-hat-expansion}
                E_0 + E_1 + E_2 + \dots + E_{m-1} + E_m \quad &\text{if $\abs{x-y} = 2m$} \\
                E_0 + E_1 + E_2 + \dots + E_{m-1} + E_m + E_{2m+1} \quad &\text{if $\abs{x-y} = 2m+1$} \,,
            \end{cases}
   \end{equation}
where $\set{E_k}_{k\geq 0}$ form a set of independent Exponential random variables, with $E_k$ having rate $2k$. (Note that $E_0\equiv 0$: it is included merely for notational convenience.)

Now define
    \begin{equation}
        \hat{v}(x,y,t) = \prob{\hat{\tau}>t \,|\, X_0 = x, \, Y_0=y}
    \end{equation}
to be the tail probability of the coupling time under $\hat{c}$. The main result of this paper is the following.
    \begin{thm}\label{thm:stoch-min}
        For any states $x,y\in\hyp$ and time $t\geq 0$,
            \begin{equation}\label{eqn:thm-statement}
                \hat{v}(x,y,t) = \inf_{c\in\mathcal{C}} \prob{\tau^c >t \,|\, X_0=x, Y_0=y} \,.
            \end{equation}
        In other words, $\hat{\tau}$ is the stochastic minimum of all co-adapted coupling times for the pair $(X,Y)$.
    \end{thm}

It is clear from the representation in \eqref{eqn:tau-hat-expansion} that $\hat{v}(x,y,t)$ only depends on $(x,y)$ through $\abs{x-y}$, and so we shall usually simply write
    \[ \hat{v}(k,t) = \prob{\hat{\tau}>t \,|\, N_0 = k} \,, \]
with the convention that $\hat{v}(k,t)=0$ for $k\leq 0$. Note, again from \eqref{eqn:tau-hat-expansion}, that $\hat{v}(k,t)$ is strictly increasing in $k$. For a strategy $c\in\mathcal{C}$, define the process $S_t^c$ by
    \[ S_t^c = \hat{v}\bra{X_t^c, Y_t^c, T-t}\,, \]
where $T>0$ is some fixed time. This is the conditional probability of $X$ and $Y$ not having coupled by time $T$, when strategy $c$ has been followed over the interval $[0,t]$ and $\hat{c}$ has then been used from time $t$ onwards. The optimality of $\hat{c}$ will follow by Bellman's principle (see, for example, \cite{Krylov-1980}) if it can be shown that $S^c_{t\wedge\tau^c}$ is a submartingale for all $c\in\mathcal{C}$, as demonstrated in the following lemma. (Here and throughout, $s\wedge t=\min\set{s,t}$.)

    \begin{lem}\label{lem:Bellman}
        Suppose that for each $c\in\mathcal{C}$ and each $T\in\RR_+$,
            \[ \bra{S_{t\wedge\tau^c}^c}_{0\leq t\leq T} \quad\text{is a submartingale.} \]
        Then equation~\eqref{eqn:thm-statement} holds.
    \end{lem}

    \begin{proof}
        Notice that, with $(X_0,Y_0)=(x,y)$, $S_0^c = \hat{v}(x,y,T)$ and $S_{T\wedge\tau^c}^c = \indev{T<\tau^c}$. If $S^c_{\cdot\,\wedge\tau^c}$ is a submartingale it follows by the Optional Sampling Theorem that
            \[ \prob{\tau^c>T} = \Ex{S^c_{T\wedge\tau^c}} \geq S_0^c = \hat{v}(x,y,T) = \prob{\hat{\tau}>T} \,, \]
        and hence the infimum in \eqref{eqn:thm-statement} is attained by $\hat{c}$.
    \end{proof}

Now, (point process) stochastic calculus yields:
    \begin{equation}\label{eqn:diff-eqn}
        dS_t^c = dZ^c_t  + \bra{\mathcal{A}_t^c \hat{v} - \frac{\partial \hat{v}}{\partial t}} dt \,,
    \end{equation}
where $Z^c_t$ is a martingale, and $\mathcal{A}_t^c$ is the {\lq\lq}generator{\rq\rq} corresponding to the matrix $Q^c(t)$. Since the Poisson processes $\Lambda_{ij}$ are independent, the probability of two or more jumps occurring in the superimposed process $\bigcup\Lambda_{ij}$ in a time interval of length $\delta$ is
$O(\delta^2)$. Hence, for any function $f:\hyp\times\hyp\times\RR^+\rightarrow\RR$, $\mathcal{A}_t^c$ satisfies
    \[ \mathcal{A}_t^c f (x,y,t) = \sum_{i=0}^n \sum_{j=0}^n q_{ij}^c(t) {\Big [} f(x+e_i,y+e_j,t) - f(x,y,t) {\Big ]} \,. \]
Setting $f=\hat{v}$ gives:
    \begin{align*}
        \mathcal{A}_t^c \hat{v}(x,y,t) &= \sum_{i=0}^n \sum_{j=0}^n q_{ij}^c(t) {\Big [} \hat{v}(x+e_i,y+
                e_j,t) - \hat{v}(x,y,t){\Big ]} \\
            &= \sum_{i=0}^n \sum_{j=0}^n q_{ij}^c(t) {\Big [}\hat{v}(\abs{x-y+e_i+
                e_j},t) - \hat{v}(\abs{x-y},t){\Big ]} \,.
    \end{align*}
In particular, since $\hat{v}$ is invariant under coordinate permutation, if $N_t^c=\abs{x-y}=k$ then    \begin{equation}\label{eqn:A-for-v-hat}
        \mathcal{A}_t^c \hat{v}(x,y,t)  = \sum_{m=-2}^2 \lambda^c_t(k,k+m) {\Big [}\hat{v}(k+m,t) -
            \hat{v}(k,t){\Big ]} \,,
    \end{equation}
where $\lambda^c_t(k,k+m)$ is the rate (according to $Q^c(t)$) at which $N_t^c$ jumps from $k$ to $k+m$. More explicitly,
    \begin{alignat}{2}
        \lambda_t^c(k,k+2) &= \sum_{\stackrel{i,j\in M_t}{i\neq j}} q_{ij}^c(t)\,, \qquad
        &\lambda_t^c(k,k+1) &= \sum_{i\in M_t} \bra{q_{i0}^c(t) + q_{0i}^c(t)}\,, \label{eqn:lambda1}\\
        \lambda_t^c(k,k-2) &= \sum_{\stackrel{i,j\in U_t}{i\neq j}} q_{ij}^c(t)\,, \qquad
        &\lambda_t^c(k,k-1) &= \sum_{i\in U_t} \bra{q_{i0}^c(t) + q_{0i}^c(t)} \,,\label{eqn:lambda2}
    \end{alignat}
and
    \begin{equation}\label{eqn:lambda3}
        \lambda_t^c(k,k) = \sum_{i\in U_t, j\in M_t} \bra{q_{ij}^c(t)+q_{ji}^c(t)} + \sum_{i=1}^n q_{ii}^c(t)\,.
    \end{equation}
It follows from the definition of $Q$ and equations~\eqref{eqn:lambda1} to \eqref{eqn:lambda3} that these terms must satisfy the linear constraints:
    \begin{align*}
        &\lambda_t^c(k,k-2) + \frac{1}{2}\lambda_t^c(k,k-1) \leq k \,, \quad\text{and} \\
        &\lambda_t^c(k,k-2) + \frac{1}{2}\lambda_t^c(k,k-1) + \lambda_t^c(k,k) + \frac{1}{2}\lambda_t^c(k,k+1) + \lambda_t^c(k,k+2) = n \,.
    \end{align*}
Denote by $L_n$ the set of non-negative $\lambda$ satisfying the constraints
    \begin{align}
        &\lambda(k,k-2) + \frac{1}{2}\lambda(k,k-1) \leq k \,, \quad\text{and} \label{eqn:lambda-constraint1} \\
        &\lambda(k,k-2) + \frac{1}{2}\lambda(k,k-1) + \lambda(k,k) + \frac{1}{2}\lambda(k,k+1) + \lambda(k,k+2) = n \,. \label{eqn:lambda-constraint2}
    \end{align}


Returning to equation~\eqref{eqn:diff-eqn}:
    \begin{equation*}
        dS_t^c = dZ^c_t  + \bra{\mathcal{A}_t^c \hat{v} - \frac{\partial \hat{v}}{\partial t}} dt \,.
    \end{equation*}
We wish to show that $S_{t\wedge\tau^c}^c$ is a submartingale for all couplings $c\in\mathcal{C}$. We shall do this by showing that $\mathcal{A}_t^c\hat{v}$ is minimised by setting $c=\hat{c}$. This is sufficient because $S^{\hat{c}}_{t\wedge\hat{\tau}}$ is a martingale (and so $\mathcal{A}_t^{\hat{c}} \hat{v} - \partial \hat{v}/\partial t =0$). Now, from equation~\eqref{eqn:A-for-v-hat} we know that
    \[ \mathcal{A}_t^c \hat{v}(k,t) = \sum_{m=-2}^2 \lambda^c_t(k,k+m) {\Big [}\hat{v}(k+m,t) - \hat{v}(k,t){\Big]
        } \,. \]
Thus we seek to show that, for all $k\geq 0$ and for all $t\geq 0$,
    \begin{equation}\label{eqn:max-aim}
        \max_{\lambda\in L_n} \sum_{m=-2}^2 \lambda(k,k+m) {\Big [}\hat{v}(k,t) - \hat{v}(k+m,t){\Big]} \geq 0 \,.
    \end{equation}
For each $t$, this is a linear function of non-negative terms of the form $\lambda(k,k+m)$. Thanks to the monotonicity in its first argument of $\hat{v}$, the terms appearing in the left-hand-side of \eqref{eqn:max-aim} are non-positive if and only if $m$ is non-negative. Hence we must set
    \begin{equation}\label{eqn:no-breaks}
        \lambda(k,k+1) = \lambda(k,k+2) = 0
    \end{equation}
in order to achieve the maximum in~\eqref{eqn:max-aim}. 

It now suffices to maximise
    \begin{equation}\label{eqn:maximisation-problem}
        \lambda(k,k-1) {\Big [}\hat{v}(k,t) - \hat{v}(k-1,t){\Big]} + \lambda(k,k-2) {\Big [}\hat{v}(k,t) - \hat{v}(k-2,t){\Big]} \,,
    \end{equation}
subject to the constraint in \eqref{eqn:lambda-constraint1}.

Combining \eqref{eqn:lambda-constraint1} and \eqref{eqn:maximisation-problem} 
yields the final version of our optimisation problem:
    \begin{align}
        \text{maximise} \qquad & \lambda(k,k-1)\bra{ {\Big [}\hat{v}(k,t) - \hat{v}(k-1,t){\Big]} -\frac{1}{2} {\Big [}\hat{v}(k,t) - \hat{v}(k-2,t){\Big]}} 
            \label{eqn:final-max}\\
        \text{subject to} \qquad & 0 \leq \lambda(k,k-1) \leq 2k \,. \label{eqn:final-constraint}
    \end{align}
The solution to this problem is clearly given by:
    \begin{equation}\label{eqn:answer}
        \lambda(k,k-1) = \begin{cases}
                                2k &\quad\text{if ${\Big [}\hat{v}(k,t) - \hat{v}(k-1,t){\Big]} > \frac{1}{2} {\Big [}\hat{v}(k,t) - \hat{v}(k-2,t){\Big]}$} \\
                                0 &\quad\text{otherwise} \,.
                            \end{cases}
    \end{equation}
These observations may be summarised as follows:

    \begin{prop}\label{prop:max-summary}
        For $\lambda\in L_n$, the maximum value of
            \[ \sum_{m=-2}^2 \lambda(k,k+m) {\Big [}\hat{v}(k,t) - \hat{v}(k+m,t){\Big]} \,, \]
        is achieved at $\lambda^*$, where $\lambda^*$ satisfies the following:
            \begin{align*}
                &\lambda^*(k,k+1)=\lambda^*(k,k+2) = 0 \,; \\
                &\lambda^*(k,k-2) + \frac{1}{2}\lambda^*(k,k-1)=k \,; \\
                &\lambda^*(k,k-1) = \begin{cases}
                                2k &\quad\text{if ${\Big [}\hat{v}(k,t) - \hat{v}(k-1,t){\Big]} > \frac{1}{2} {\Big [}\hat{v}(k,t) - \hat{v}(k-2,t){\Big]}$} \\
                                0 &\quad\text{otherwise} \,.
                            \end{cases}
            \end{align*}

    \end{prop}

\medskip
\noindent Our final proposition shows that $\lambda^*(k,k-1) = 2k$ if and only if $k$ is odd.
    \begin{prop}\label{prop:parity}
        For any fixed $t\geq 0$,
         \begin{align}
            2{\Big [}\hat{v}(k,t) - \hat{v}(k-1,t){\Big]} - {\Big [}\hat{v}(k,t) - \hat{v}(k-2,t){\Big]}
                    &\geq 0 \quad\text{if $k$ is odd, and} \label{eqn:prop-parity1} \\
            2{\Big [}\hat{v}(k,t) - \hat{v}(k-1,t){\Big]} - {\Big [}\hat{v}(k,t) - \hat{v}(k-2,t){\Big]}
                    &\leq 0 \quad\text{if $k$ is even.} \label{eqn:prop-parity2}
         \end{align}
    \end{prop}

    \begin{proof}
        Define $\hat{V}_\alpha$ by
            \begin{equation*}
                \hat{V}_\alpha(k) = \int_0^\infty e^{-\alpha t} \hat{v}(k,t) dt
                    = \frac{1}{\alpha}\bra{1-\Ex{e^{-\alpha\hat{\tau}}}} \,.
            \end{equation*}
        We also define $d(k,t) = \hat{v}(k,t)-\hat{v}(k-1,t)$, and for $\alpha\geq 0$ let
            \[ D_\alpha(k) = \int_0^\infty e^{-\alpha t} d(k,t) dt \]
        be the Laplace transform of $d(k,\cdot)$.
        Given the representation in equation~\eqref{eqn:tau-hat-expansion} of $\hat{\tau}$ as a sum of independent Exponential random variables, it follows that
            \begin{equation}\label{eqn:LT-for-v-hat}
                \hat{V}_\alpha(k) = \left\{
                    \begin{alignedat}{2}
                        &\frac{1}{\alpha}\bra{1- \prod_{i=1}^m \frac{2i}{2i+\alpha}} &\quad &\text{if $k=2m$}
                            \\
                        &\frac{1}{\alpha}\bra{1- \frac{2(2m+1)}{2(2m+1)+\alpha} \prod_{i=1}^m
                            \frac{2i}{2i+\alpha}} &\quad &\text{if $k=2m+1$} \,.
                    \end{alignedat}
                                    \right.
            \end{equation}
        To ease notation, let
            \[ \phi_\alpha(m) = \prod_{i=1}^m \frac{2i}{2i+\alpha} \,. \]
        The following equality then follows directly from consideration of the transition rates corresponding to strategy $\hat{c}$: \\
        \indent for all $\alpha\geq 0$ and $m\geq 1$,
            \begin{align}
                1 - \alpha \hat{V}_\alpha(2m) + 2m \squ{\hat{V}_\alpha(2m-2) - \hat{V}_\alpha(2m)}
                &= \phi_\alpha(m) + \frac{2m}{\alpha} \squ{\phi_\alpha(m)-\phi_\alpha(m-1)} \nonumber \\
                &= \phi_\alpha(m) + \frac{2m}{\alpha}\phi_\alpha(m) \squ{1-\frac{2m+\alpha}{2m}} \nonumber \\
                &= 0 \,. \label{eqn:imp-V1}
            \end{align}
        Similarly,
            \begin{align}
                1 &- \alpha \hat{V}_\alpha(2m-1) + 2(2m-1) \squ{\hat{V}_\alpha(2m-2) - \hat{V}_\alpha(2m-1)} = 0 \,. \label{eqn:imp-V2}
            \end{align}

        Now suppose that $k=2m$, and hence is even. We wish to prove that
            \[ d(2m-1,t)-d(2m,t) \geq 0 \quad\text{for all $t\geq 0$}\,, \]
        which is equivalent to showing that $D_\alpha(2m-1) - D_\alpha(2m)$ is totally (or completely) monotone (by the Bernstein-Widder Theorem; Theorem 1a of \cite{Feller-1971}, Ch. XIII.4).

        We proceed by subtracting equation~\eqref{eqn:imp-V2} from \eqref{eqn:imp-V1}:
            \begin{align*}
                    0&= -\alpha \squ{\hat{V}_\alpha(2m)-\hat{V}_\alpha(2m-1)} + 2m
                        \squ{\hat{V}_\alpha(2m-2)-\hat{V}_\alpha(2m)} \\
                    &\phantom{= -\alpha \squ{\hat{V}_\alpha(2m)-\hat{V}_\alpha(2m-1)} } + 2(2m-1)\squ{\hat{V}_\alpha(2m-1)-\hat{V}_\alpha(2m-2)} \\
                    &= -\alpha D_\alpha(2m)-2m\squ{D_\alpha(2m)+D_\alpha(2m-1)} + 2(2m-1)D_\alpha(2m-1)\,,
            \end{align*}
        and so
            \begin{equation}\label{eqn:D-even}
              D_\alpha(2m-1)-D_\alpha(2m)= \frac{2+\alpha}{2m-2}D_\alpha(2m)  \,.
            \end{equation}
        It therefore suffices to show that $(2+\alpha)D_\alpha(2m)$ is completely monotone.

        Now note from the form of $\hat{V}$ in equation~\eqref{eqn:LT-for-v-hat}, that
            \[ (2+\alpha)D_\alpha(2m) = 2\Theta_\alpha(2m) \,, \]
        where $\Theta_\alpha(2m)$ is the Laplace transform of
            \[ \theta(2m,t) = \prob{\sum_{i=0}^m E_i >t} - \prob{\sum_{i=0}^{m-1} E_i + E_{2m-1} >t} \,, \]
        where $\set{E_i}_{i\geq 0}$ form a set of independent Exponential random variables, with $E_i$ having parameter $2i$. But since $\theta(2m,t)$ is strictly positive for all $t$, it follows that \\ $(2+\alpha)D_\alpha(2m)$ is completely monotone, as required. This proves that, for any fixed $t\geq 0$,
            \begin{equation}\label{eqn:prop-converse}
                2{\Big [}\hat{v}(k,t) - \hat{v}(k-1,t){\Big]} - {\Big [}\hat{v}(k,t) - \hat{v}(k-2,t){\Big]} \leq 0
            \end{equation}
        whenever $k$ is even. Thus inequality~\eqref{eqn:prop-parity2} holds in this case.

        Now suppose that $k=2m+1$, and hence is odd. In this case we wish to show that inequality~\eqref{eqn:prop-parity1} holds, which is equivalent to showing that $D_\alpha(2m+1) - D_\alpha(2m)$ is completely monotone. Now, substituting $m+1$ for $m$ in equation~\eqref{eqn:imp-V2} yields
            \begin{equation}\label{eqn:imp-V2b}
                 1 - \alpha \hat{V}_\alpha(2m+1) + 2(2m+1) \squ{\hat{V}_\alpha(2m) - \hat{V}_\alpha(2m+1)} = 0 \,.
            \end{equation}
        Proceeding as above, we subtract equation~\eqref{eqn:imp-V1} from \eqref{eqn:imp-V2b}:
            \begin{align}
              0 &= -\alpha \squ{\hat{V}_\alpha(2m+1)-\hat{V}_\alpha(2m)} + 2(2m+1)
                    \squ{\hat{V}_\alpha(2m)-\hat{V}_\alpha(2m+1)} \nonumber \\
                &\phantom{= -\alpha \squ{\hat{V}_\alpha(2m+1)-\hat{V}_\alpha(2m)} } +2m\squ{\hat{V}_\alpha(2m)-\hat{V}_\alpha(2m-2)}\nonumber  \\
                &= -\alpha D_\alpha(2m+1)-2(2m+1)D_\alpha(2m+1) + 2m\squ{D_\alpha(2m)+D_\alpha(2m-1)} \,. \label{eqn:difference}
            \end{align}
       Then it follows from equation~\eqref{eqn:D-even} that
            \begin{equation}\label{eqn:ident}
                (2m-2)D_\alpha(2m-1) = (2m+\alpha)D_\alpha(2m) \,.
            \end{equation}
       Substitution of equation~\eqref{eqn:ident} into \eqref{eqn:difference} gives
            \[ 0 = (4m+2-\alpha)\squ{D_\alpha(2m)-D_\alpha(2m+1)} +2\squ{D_\alpha(2m-1)-D_\alpha(2m)} \,, \]
       and so
            \begin{equation}\label{eqn:D-odd}
              D_\alpha(2m+1)-D_\alpha(2m)= \frac{2}{4m+2+\alpha}\squ{D_\alpha(2m-1)-D_\alpha(2m)} \,.
            \end{equation}
       But, since we have already seen that $D_\alpha(2m-1)-D_\alpha(2m)$ is completely monotone, the right-hand-side of equation~\eqref{eqn:D-odd} is the product of two completely monotone functions, and so is itself completely monotone~\cite{Feller-1971}, as required.
    \end{proof}

Now we may complete the
	\begin{proof}[Proof of Theorem~\ref{thm:stoch-min}]
Thanks to Lemma~\ref{lem:Bellman} and Proposition~\ref{prop:max-summary}, Proposition~\ref{prop:parity}, along with equations~\eqref{eqn:no-breaks} and \eqref{eqn:answer}, shows that any optimal choice of $Q(t)$, $Q^*(t)$, is of the following form:
    \begin{itemize}
        \item when $N_t$ is odd:
            \begin{align*}
                 q^*_{i0}(t) &= q^*_{0i}(t) = 1 \,\,\text{for all $i\in U_t$, (and so $\lambda_t^*(N_t,N_t-1)=2N_t$)} \,, \\
                 q^*_{ii}(t) &= 1 \,\,\text{for all $i\in M_t$} \,;
            \end{align*}
        \item when $N_t$ is even:
            \begin{align}
                 q^*_{i0}(t) &= q^*_{0i}(t) =q^*_{ii}(t) = 0 \,\,\text{for all $i\in U_t$, (and so $\lambda_t^*(N_t,N_t-1)=0$)} \,, \label{eqn:Q*} \\
                 q^*_{ii}(t) &= 1 \,\,\text{for all $i\in M_t$} \nonumber \,.
            \end{align}
    \end{itemize}
This is in agreement with our candidate strategy $\hat{Q}$ (recall Definition~\ref{defn:Q-hat}). From equation~\eqref{eqn:Q*} it follows that the values of $q^*_{ij}(t)$ for distinct $i,j\in U_t$ must satisfy
    \[ \sum_{\stackrel{i,j\in U_t}{i\neq j}} q^*_{ij}(t) = \abs{U_t} \,, \]
but are not constrained beyond this. Our choice of
    \[ \hat{q}_{ij}(t) = \frac{1}{\abs{U_t}-1} \]
satisfies this bound, and so $\hat{c}$ is truly an optimal co-adapted coupling, as claimed.
	\end{proof}

	\begin{rem}
Observe that when $k=1$, equation~\eqref{eqn:tau-hat-expansion} implies that $\hat{v}(1,t)=\hat{v}(2,t)$ for all $t$. The optimisation problem in \eqref{eqn:final-max} and \eqref{eqn:final-constraint} simplifies in this case to the following:
    \begin{align}
        \text{maximise} \qquad & \lambda(1,0)\hat{v}(1,t) \label{eqn:final-max2}\\
        \text{subject to} \qquad & \frac{1}{2}\lambda(1,0) + \lambda(1,1) +\frac{1}{2}\lambda(1,2)\leq n\, .\label{eqn:final-constraint2}
    \end{align}
As above, this is achieved by setting $\lambda(1,0)=2$. Note from equation~\eqref{eqn:final-constraint2}, however, that when $k=1$ there is no obligation to set $\lambda(1,2)=0$ in order to attain the required maximum. Indeed, due to the equality between $\hat{v}(1,t)$ and $\hat{v}(2,t)$, when $k=1$ it is not sub-optimal to allow \emph{matched} coordinates to evolve independently (corresponding to $\lambda_t^c(1,2)>0$), so long as strategy $\hat{c}$ is used once more as soon as $k=2$.
	\end{rem}


\section{Maximal coupling}\label{sec:maximal}

Let $X$ and $Y$ be two copies of a Markov chain on a countable space, starting from different states. The coupling inequality (see, for example, \cite{Lindvall-2002}) bounds the tail distribution of \emph{any} coupling of $X$ and $Y$ by the total variation distance between the two processes:
    \begin{equation}\label{eqn:coupling-inequality}
        \norm{\law(X_t)-\law(Y_t)}_{TV} \leq \prob{\tau>t} \,.
    \end{equation}
Griffeath \cite{Griffeath-1975} showed that, for discrete-time chains, there always exists a \emph{maximal} coupling of $X$ and $Y$: that is, one which achieves equality for all $t\geq 0$ in the coupling inequality. This result was extended to general continuous-time stochastic processes with paths in Skorohod space in \cite{Sverchkov.Smirnov-1990}. However, in general such a coupling is not co-adapted. In light of the results of Section~\ref{sec:optimal}, where it was shown that $\hat{c}$ is the optimal co-adapted coupling for the symmetric random walk on \hyp, a natural question is whether $\hat{c}$ is also a maximal coupling.

This is certainly not the case in general. Suppose that $X$ and $Y$ are once again random walks on \hyp, with $X_0=(0,0,\dots,0)$ and $Y_0=(1,1,\dots,1)$: calculations as in \cite{Diaconis.Graham.ea-1990} show that the total variation distance between $X_t$ and $Y_t$ exhibits a cutoff phenomenon, with the cutoff taking place at time $T_n=\frac{1}{4}\log n$ for large $n$. This implies that a maximal coupling of $X$ and $Y$ has expected coupling time of order $T_n$. However, it follows from the representation of $\hat{\tau}$ in equation~\eqref{eqn:tau-hat-expansion} that
    \begin{equation}\label{eqn:expectation-of-tau}
        \Ex{\,\hat{\tau}\,\,;\,\abs{X_0-Y_0}=n=2m} = \Ex{E_1 + E_2 + \dots + E_{m-1} + E_m} \sim \frac{1}{2}\log(n) \,.
    \end{equation}
It follows that $\hat{c}$ is not, in general, a maximal coupling.

A faster coupling of $X$ and $Y$ was proposed by \cite{Matthews-1987}. This coupling also makes new coordinate matches in pairs, but uses information about the future evolution of one of the chains in order to make such matches in a more efficient manner. This coupling is very near to being maximal (it captures the correct cutoff time), but is of course not co-adapted. Further results related to the construction of maximal couplings for general Markov chains may be found in \cite{Greven-1987,Harison.Smirnov-1990,Pitman-1976}.

\bibliographystyle{apt}
\bibliography{AP-12589-Connor-Jacka}

\begin{thebibliography}{10}

\bibitem{Aldous-1983}
{\sc Aldous, D.} (1983).
\newblock {\em Random walks on finite groups and rapidly mixing {M}arkov
  chains} vol.~986 of {\em Lecture Notes in Math.}
\newblock Springer, Berlin.
\newblock pp.~243--297.

\bibitem{Diaconis.Graham.ea-1990}
{\sc Diaconis, P., Graham, R.~L. and Morrison, J.~A.} (1990).
\newblock Asymptotic analysis of a random walk on a hypercube with many
  dimensions.
\newblock {\em Random Structures Algorithms\/} {\bf 1,} 51--72.

\bibitem{Feller-1971}
{\sc Feller, V.} (1971).
\newblock {\em An introduction to probability theory and its applications}
  second~ed. vol.~2.
\newblock Wiley.

\bibitem{Greven-1987}
{\sc Greven, A.} (1987).
\newblock Couplings of {M}arkov chains by randomized stopping times.
\newblock {\em Probab. Theory Related Fields\/} {\bf 75,} 195--212 and
  431--458.

\bibitem{Griffeath-1975}
{\sc Griffeath, D.} (1975).
\newblock A maximal coupling for {M}arkov chains.
\newblock {\em Z. Wahrscheinlichkeitstheorie und Verw. Gebiete\/} {\bf 31,}
  95--106.

\bibitem{Harison.Smirnov-1990}
{\sc Harison, V. and Smirnov, S.~N.} (1990).
\newblock Jonction maximale en distribution dans le cas markovien.
\newblock {\em Probab. Theory Related Fields\/} {\bf 84,} 491--503.

\bibitem{Krylov-1980}
{\sc Krylov, N.~V.} (1980).
\newblock {\em Controlled diffusion processes} vol.~14 of {\em Applications of
  Mathematics}.
\newblock Springer-Verlag, New York.
\newblock Translated from the Russian by A. B. Aries.

\bibitem{Lindvall-2002}
{\sc Lindvall, T.} (2002).
\newblock {\em Lectures on the coupling method}.
\newblock Dover.

\bibitem{Matthews-1987}
{\sc Matthews, P.} (1987).
\newblock Mixing rates for a random walk on the cube.
\newblock {\em SIAM J. Algebraic Discrete Methods\/} {\bf 8,} 746--752.

\bibitem{Pitman-1976}
{\sc Pitman, J.~W.} (1976).
\newblock On coupling of {M}arkov chains.
\newblock {\em Z. Wahrscheinlichkeitstheorie und Verw. Gebiete\/} {\bf 35,}
  315--322.

\bibitem{Sverchkov.Smirnov-1990}
{\sc Sverchkov, M.~Y. and Smirnov, S.~N.} (1990).
\newblock Maximal coupling of {D}-valued processes.
\newblock {\em Soviet Math. Dokl.\/} {\bf 41,} 352--354.

\end{thebibliography}

\end{document}